\input amstex
\documentstyle{amsppt}
\document
\topmatter
\title
Higher dimensional Gray Hermitian manifolds
\endtitle
\author
W\l odzimierz Jelonek
\endauthor

\abstract{The aim of this paper is to describe a large class of
Hermitian Gray manifolds. }
\thanks{MS Classification: 53C55, 53C25. Key words and phrases:
Gray manifold, Hermitian metric, Killing tensor, special
K\"ahler-Ricci potential }\endthanks
 \endabstract
\endtopmatter
\define\G{\Gamma}
\define\DE{\Cal D^{\perp}}
\define\e{\epsilon}
\define\n{\nabla}

\define\w{\wedge}
\define\k{\diamondsuit}
\define\th{\theta}
\define\p{\partial}
\define\a{\alpha}

\define\g{\gamma}
\define\lb{\lambda}

\define\AC{\Cal AC^{\perp}}

\define\1{D_{\lb}}
\define\2{D_{\mu}}
\define\0{\Omega}

\define\De{\Cal D}

\define\m{(M,g,J)}
\define \E{\Cal E}
\bigskip
{\bf 0. Introduction. }  The Gray manifolds (also called
$\AC$-manifolds)
 are Riemannian manifolds $(M,g)$ whose
 Ricci tensor $\rho$ satisfies the condition that $\rho-\frac{2\tau}{n+2}g$
is a Killing tensor, where $\tau$ is the scalar curvature of
$(M,g)$. The equivalent condition is
$\n_X\rho(X,X)-\frac{2X\tau}{n+2}g(X,X)=0$ for all $X\in TM$. This
condition was first considered by A. Gray in [G]. The class of
Gray manifolds is in some sense close to Einstein manifolds.
 The problem of
 describing compact Gray manifolds was  stated by Besse in [Bes].
 The Einstein structures on $S^2$-bundles were constructed by
J. Wang and M. Wang in [W-W] and B\'erard Bergery in [Ber],
non-homogeneous K\"ahler Einstein manifolds were constructed in
[K-S], and conformally K\"ahler Einstein manifolds were classified
by Derdzi\'nski and Maschler in [D-M-2].  In  [J-2],  [J-3] we
have described the class of compact four
 dimensional Gray manifolds  with two eigenvalues of the
 Ricci tensor.  These manifolds were the sphere bundles over
 Riemannian surfaces of any genus $g$. In addition  these
 manifolds were  Hermitian with respect to two opposite complex
  structures.
 The Gray K\"ahler manifolds are weakly Bochner flat -(WBF) K\"ahler manifolds.
 The class of WBF K\"ahler manifolds was investigated for example in ([A-C-G], [A-C-G-T]).

 The aim of the present paper is to classify compact Gray
 manifolds with two eigenvalues of the Ricci tensor,  whose Ricci
 tensor is Hermitian with respect to a conformally K\"ahler complex
 structure and which have two-dimensional, integrable and totally
 geodesic eigendistribution of the  Ricci tensor.  Using the
 theorem of Derdzi\'nski and Maschler ([D-M-1], [D-M-2]) we prove that such
 structures exist  on $\Bbb{CP}^1$-bundles over K\"ahler-Einstein manifolds of any
dimension $2(n-1)$, $n\ge2$, and on $\Bbb{CP}^{n}$ and give a
complete classification of such metrics.
  We shall prove that in the case of $\Bbb{CP}^1$-bundles the eigenvalues of
the Ricci tensor are everywhere different and our manifolds are
also  Hermitian with respect to an  opposite oriented complex
structure (we call them  bi-Hermitian). In the case of
$\Bbb{CP}^n$ the eigenvalues coincide  at exactly one point where
the opposite Hermitian structure is not defined. In particular, we
classify  the class of compact Gray K\"ahler manifolds (i.e. WBF
K\"ahler manifolds) whose Ricci tensor has two eigenvalues and
two-dimensional, totally geodesic eigendistribution. These
K\"ahler manifolds are also described in [A-C-G-T].

 \par
\bigskip
{\bf 1. Killing tensors.}  Our notation is as in [K-N]. Let
$(M,g)$ be a smooth connected Riemannian manifold. For a tensor
$T(X_1,X_2,\ldots,X_k)$ we define a tensor $\n T(X_0,X_1,
\ldots,X_k)$ by $\n T(X_0,X_1, \ldots,X_k)=\n_{X_0}T(X_1,
\ldots,X_k)$.
 By a Killing tensor on $M$
we mean  an endomorphism $S\in$ End$(TM)$ satisfying the following
conditions:

(1.1)  $ g(SX,Y)=g(X,SY)$  for all $X,Y\in TM$,

(1.2)  $g(\n S(X,X),X)=0$ for all $X \in TM$.

 We call
$S$
 a proper Killing tensor if $\n S\ne 0$.
We  denote by $\Phi$ the tensor defined by $\Phi(X,Y)=g(SX,Y).$

We start with :
\par
\medskip
  {\bf Proposition 1.1.} {\it The following conditions are equivalent:}

  (a) {\it A tensor} $S$ {\it is a Killing tensor on} $(M,g)$;

  (b) {\it For every geodesic }$\gamma$ {\it on} $(M,g)$ {\it
    the function}
  $\Phi (\g '(t),\g '(t))$ {\it is constant on }  dom$\g$.

  (c)  {\it The  condition}
  $$ \n_X \Phi(Y,Z)+ \n_Z \Phi(X,Y)+\n_Y \Phi(Z,X)=0\tag A$$
  {\it is satisfied for all }$X,Y,Z\in \frak X(M)$.
\par
\medskip

  {\it Proof:} By using polarization it is easy to see that (a) is
  equivalent to (c). Let now  $X \in T_{x_0}M$ be any vector
  and $\g$ be a geodesic satisfying the initial condition $\g '(0)=X$.
  Then
  $$ \frac d{dt}\Phi(\g '(t),\g '(t))=\n_{\g '(t)}\Phi(\g '(t),\g '(t)).
  \tag 1.3$$
  Hence $ \frac d{dt}\Phi(\g '(t),\g '(t))_{t=0}=\n\Phi(X,X,X).$
  The equivalence (a) $\Leftrightarrow$ (b) follows immediately
  from the above relations.$\k$
  \par
  \medskip
   As in [D] define   $E_S(x)$ to be the number
   of distinct eigenvalues of $S_x$, and set $M_S=\{ x\in M:E_S$ is
   constant in a neighbourhood of $x\}$. The set $M_S$ is open and dense
   in $M$ and the eigenvalues $\lb_i$ of $S$ are distinct and smooth
   in each component $U$ of $M_S$. The eigenspaces $D_{\lb}=\ker (S-\lb I)$
    form smooth distributions in each component $U$ of $M_S$. By
    $\n f$ we denote the gradient of a function $f$ ($g(\n f,X)=df(X)$)
    and by $\G(D_{\lb})$  (resp.  $\frak X(U)$) the set of all local
    sections of the bundle
    $D_{\lb}$ (resp. of all local vector fields on $U$). Let us note that if $\lb\ne \mu$ are eigenvalues
    of $S$ then $\1$ is orthogonal to $\2$.
    \par
    \medskip

   {\bf Theorem 1.2.} {\it Let } $S$ {\it be a Killing  tensor on} $M$
   {\it and U be a component of } $M_S$ {\it and}
   $\lb_1,\lb_2,..,\lb_k \in C^{\infty}(U)$ {\it be eigenfunctions of }
   $S$. {\it Then for all} $X\in D_{\lb_i}$ {\it we have}
   $$ \n S(X,X)=-\frac12 \n \lb_i\| X \| ^2\tag 1.4$$
   {\it and}  $D_{\lb_i}\subset \ker d\lb_i$. { \it If } $i\ne j$
   {\it and} $X\in \G(D_{\lb_i}),Y \in \G(D_{\lb_j})$
   {\it then}

  $$ g(\n_X X ,Y)=\frac12\frac{Y\lb_i}{\lb_j-\lb_i}\| X \| ^2.
   \tag 1.5$$
\par
\medskip
   {\it Proof.} Let $X \in \G(D_{\lb_i})$ and $Y \in \frak X(U).$
   Then we have $SX= \lb_iX$ and
   $$  \n S(Y,X)+(S-\lb_iI)(\n_YX)=(Y\lb_i)X \tag 1.6$$
   and consequently,
   $$  g(\n S(Y,X),X)=(Y\lb_i)\| X\|^2. \tag 1.7$$
   Taking $Y=X$ in (1.7) we obtain  $0=X\lb_i\| X\|^2$ by (1.2).
   Hence $D_{\lb_i}\subset \ker d\lb_i$. Thus from (1.6) it follows that
   $ \n S(X,X)=(\lb_iI-S)(\n_XX)$. Condition  (A) implies
   $ g(\n S(X,Y),Z)+ g(\n S(Z,X),Y)+g(\n S(Y,Z),X)=0$
   hence,
    $$ 2g(\n S(X,X),Y)+ g(\n S(Y,X),X)=0.\tag 1.8 $$
    Thus, (1.8) yields
    $ Y\lb_i\parallel X\parallel^2 + 2g(\n S(X,X),Y)=0$. Consequently,
    $\n S(X,X)=-\frac12 \n \lb_i\parallel X\parallel^2$.
     Let now $Y\in \G(D_{\lb_j}).$
    Then we have
    $$  \n S(X,Y)+(S-\lb_jI)(\n_XY)=(X\lb_j)Y. \tag 1.9$$
    It is also clear that $g(\n S(X,X),Y)=g(\n S(X,Y),X)=
    (\lb_j-\lb_i)g(\n_XY,X)$. Thus,
    $$ Y\lb_i\parallel X \parallel^2=-2(\lb_j-\lb_i)g(\n_XY,X)
    =2(\lb_j-\lb_i)g( Y,\n_XX)$$
    and (1.5) holds.$\k$

    \par
    \medskip
    {\it Corrolary 1.3.} {\it Let $S,U,\lb_1,\lb_2,..,\lb_k$ be as above and
    $i\in \{1,2,..,k\}$. Then
    the following conditions are equivalent:

    (a)  For all $X\in \G(D_{\lb_i}), \ \ \n_XX \in D_{\lb_i}.$

    (b)  For all $X,Y\in \G(D_{\lb_i}), \ \ \n_XY+\n_YX \in D_{\lb_i}.$

    (c)  For all $X\in \G(D_{\lb_i}), \ \ \n S(X,X) =0.$

    (d)  For all $X,Y\in \G(D_{\lb_i}), \ \ \n S(X,Y)+\n S(Y,X) =0.$

    (e)  $\lb_i$ is a constant eigenvalue of $S$.}
    \par
    \medskip
    Let us note that if $X,Y\in \G(D_{\lb_i})$ then
    ($D_{\lb_i}\subset\ker d\lb_i $ !)
    $$\n S(X,Y)-\n S(Y,X)=(\lb_iI-S)([X,Y]), \tag 1.10$$
    hence the distribution $D_{\lb_i}$ is integrable if and only if
    $\n S(X,Y)$ $=\n S(Y,X)$ for all $X,Y\in \G(D_{\lb_i}).$  Consequently,
    we obtain
    \par
    \medskip
    {\it Corollary 1.4.} {\it Let $\lb_i \in C^{\infty}(U)$ be an
    eigenvalue of an $\Cal A$-tensor $S$. Then on $U$ the following conditions are
    equivalent:

    (a) $D_{\lb_i}$ is integrable and $\lb_i$ is constant.

    (b) For all $X,Y\in \G(D_{\lb_i}),  \ \ \n S(X,Y)=0.$

    (c) $D_{\lb_i}$ is autoparallel.}
    \par
    \medskip
    {\it Proof.}  This follows from (1.4), (1.10), Corollary 1.3 and the
    relation $\n_XY=\n_YX+[X,Y].$ $\k$

    \par
    \medskip
    Now we shall characterize  Killing tensors with two eigenvalues.
We start with:
\par
\medskip
{\bf Lemma.} {\it Let S be a self-adjoint tensor on }$(M,g)$ {\it
with exactly two eigenvalues} $\lb,\mu$. {\it If the distributions
} $D_{\lb},\ D_{\mu}$ {\it are both umbilical, $\n\lb\in\G(\2),\n
\mu\in\G(\1)$ and the mean curvatures $\xi_{\lb},\xi_{\mu}$ of the
distributions $\1,\2$ respectively satisfy the equations}
$$\xi_{\lb}=\frac1{2(\mu-\lb)}\n\lb,\qquad\xi_{\mu}=\frac1{2(\lb-\mu)}\n\mu,$$
{\it then $S$ is a Killing tensor.}
\medskip
By $\rho$ we shall denote the Ricci tensor of $(M,g)$ and by
$\tau=\operatorname{tr}_g\rho$ the scalar curvature of $(M,g)$.
\medskip
{\it Definition.}  A Riemannian manifold $(M,g)$ will be  called a
Gray $\Cal {AC}^{\perp}$ manifold if the tensor
$\rho-\frac{2\tau}{n+2}g$ is a Killing tensor.

In this paper Gray $\Cal {AC}^{\perp}$ manifolds will be  called
for short Gray manifolds or $\Cal{AC}^{\perp}$ manifolds.
\medskip
{\bf Proposition 1. } {\it Let $(M,g)$ be a $2n$-dimensional
Riemannian manifold whose Ricci tensor $\rho$ has two eigenvalues
$\lb_0(x),\mu_0(x)$ of   multiplicity 2 and $2(n-1)$ respectively
at every point $x$ of $M$. Assume that the eigendistribution
$\De_{\lb}=\De$ corresponding to $\lb_0$ is a totally geodesic
foliation and the eigendistribution $\De_{\mu}=\De^{\perp}$
corresponding to $\mu_0$ is umbilical. Then  $(M,g)$ is a Gray
manifold if and only if $\lb_0-2\mu_0$ is constant and
$\n\tau\in\G(\De)$. }
\bigskip
{\it Proof. }  Let $S_0$ be the Ricci endomorphism of $(M,g)$,
i.e. $\rho(X,Y)=g(S_0X,Y)$. Let  $S$ be the tensor defined by the
formula
$$S_0=S+\frac{\tau}{n+1} \operatorname{id}.\tag 1.11$$
Then $$\operatorname{tr} S=-\frac{(n-1)\tau}{n+1}.\tag 1.12$$ Let
$\lb_0,\mu_0$ be the eigenfunctions of $S_0$ and let us assume
that
$$\lb_0-2\mu_0=\frac{n+1}{n-1}C\tag 1.13$$
where $C\in\Bbb{R}$. Note that $S$ also has two eigenfunctions
which we denote by $\lb,\mu$ respectively. It is easy to see that
$\lb=C,\mu=-\frac{\tau}{2(n+1)}-\frac{C}{n-1}$ and
$\lb_0=\frac{\tau}{n+1}+C,\mu_0=\frac{\tau}{2(n+1)}-\frac{C}{n-1}$.
Since the distribution $\DE$ is umbilical we have
$\n_XX_{|\De}=g(X,X)\xi$ for any $X\in\G(\DE)$ where $\xi$ is the
mean curvature normal of $\DE$. Since the distribution $\De$ is
totally geodesic we also have $\n_XX_{|\DE}=0$ for any
$X\in\G(\De)$. Let $\{E_1,E_2,E_3,E_4,...,E_{2n-1},E_{2n}\}$ be a
local  orthonormal basis of $TM$ such that $\De=\text{span
}\{E_1,E_2\}$ and $\DE=\text{span }\{E_3,E_4,...,E_{2n}\}$. Then
$\n_{E_i}E_{i|\De}=\xi$ for $i\in\{3,4,...,2n\}$. Consequently
(note that $\n\mu_{|\DE}=0$ if and only if $\n\tau_{|\DE}=0$),
$$\gather \operatorname{tr}_g\n S=\sum_{i=3}^n\n S(E_i,E_i)=-2(n-1)
(S-\mu \operatorname{id})(\n_{E_3}E_{3})+\n\mu_{|\DE}\tag 1.14\\=-2(n-1)(\lb-\mu)\xi\endgather$$
if we assume that $\n\tau_{|\DE}=0$. On the other hand,
$\operatorname{tr}_g\n S_0=\frac{\n\tau}2$ and
$\operatorname{tr}_g\n S=\operatorname{tr}_g\n
S_0-\frac{\n\tau}{n+1}$. Consequently,
$$\operatorname{tr}_g\n S=\frac{(n-1)\n\tau}{2(n+1)}=-(n-1)\n\mu.\tag 1.15$$
Thus $\xi=-\frac1{2(\mu-\lb)}\n\mu$. From  the Lemma it follows
that $(M,g)$ is an $\AC$-manifold if $\lb_0-2\mu_0$ is constant
and $\n\tau\in\G(\De)$. These conditions are also necessary since
$\n\lb=0$ if $(M,g)$ is an $\AC$-manifold and $\1$ is totally
geodesic. Analogously $\xi=-\frac1{2(\mu-\lb)}\n\mu$ and
$\n\mu=-\frac1{2(n+1)}\n\tau\in\G(\De)$, where $\xi$ is the mean
curvature normal of the umbilical distribution $\2$, if $(M,g)$ is
an $\AC$-manifold.$\k$

\bigskip {\bf 2. $\AC$ Hermitian manifolds with Hermitian Ricci
tensor and two eigenvalues of the Ricci tensor.} In this section
we classify compact $\AC$ manifolds $(M,g,J)$ with two eigenvalues
of the Ricci tensor whose Ricci tensor is Hermitian with respect
to a conformally K\"ahler complex structure $J$. We shall
additionally assume that $\dim\1= 2$ in the set $U=\{x\in
M:E(x)=2\}$ where $\1$ is defined, and that the distribution $\1$
is integrable and totally geodesic. This is equivalent to the fact
that the tensor $S$ has a constant eigenvalue $\lb$ and $\1$ is
integrable. Note that the set $U$ is open and  denote by $V$ the
interior of $\{x:E(x)=1\}$. Then $U\cup V$ is an open dense subset
of $M$. It is also clear that $(V,g)$ is an Einstein manifold.  We
shall later prove that in fact $V=\emptyset$ and either $U=M$ or
$M=\Bbb{CP}^n$.

Let $f\in C^{\infty}(M)$ be a positive function such that the
manifold $(M,f^2g,J)$ is a  K\"ahler manifold.  This means that
the metric $g_1=f^2g$ is K\"ahlerian with respect to $J$.  By
$\n,\n^1$ we shall denote the Levi-Civita connections of the
metrics $g,g_1$ respectively, and by $\n F,\n^1 F$  the gradients
of a function $F$ with respect to the metrics $g,g_1$. Note that
$$\n^1F=f^{-2}\n F,\tag 2.1$$
and also
$$\n_XY=\n^1_XY-d\ln f(X)Y-d\ln f(Y)X+g_1(X,Y)\n^1\ln f.\tag 2.2$$

Note that $\n_X=\n^1_X+K_X$ where $K_X(Y)=-d\ln f(X)Y-d\ln
f(Y)X+g_1(X,Y)\n^1\ln f$. Now we compute the covariant derivative
$\n J$ of the complex structure $J$:
$$\gather  \n_XJ(Y)=[K_X,J]Y=-d\ln f(JY)X+d\ln f(Y)JX\tag 2.3\\+g_1(X,JY)\n^1\ln
f-g_1(X,Y)J\n^1\ln f.\endgather$$

Note that $g_1(X,Y)\n^1\ln f=g(X,Y)\n\ln f$.   From (2.3) it is
clear that
$$\gather  \n_XJ(X)=-d\ln f(JX)X+d\ln f(X)JX-g_1(X,X)J\n^1\ln f.\tag 2.4\endgather$$

In particular, $$\operatorname{tr}_gJ\n_XJ(X)=2(n-1)\n\ln f.\tag
2.5$$

Let $\rho,\rho_1$ be the Ricci tensors of $(M,g),(M,g_1)$
respectively. Then
$$\rho=\rho_1+(2n-2)\frac1fH^f-\bigg(\frac1f\Delta
f+(2(n-1)\frac1{f^2}\|\n^1f\|^2\bigg)g_1,\tag2.6$$ where
$H^f(X,Y)=g_1(\n^1_X\n^1f,Y)$ is the Hessian of $F$ with respect
to the metric $g_1$. Since both $\rho,\rho_1$ are $J$-invariant it
follows that the field $\xi=J(\n^1 f)$ is a holomorphic Killing
vector field for  both $(M,g,J)$ and $(M,g_1,J)$  (see [D]). Now
we prove that $\xi\in \G(\1)$. Since $\1$ is totally geodesic it
follows that for any $X,Y\in\1$ we have $\n S(X,Y)=0$. On the
other hand, $S\circ J=J\circ S$ and
$$\n_XJ\circ S+J\circ\n_XS=\n_XS\circ J+S\circ\n_XJ.\tag 2.7$$
Consequently, for $X\in \1$ we get
$$\lb\n J(X,X)+J\circ\n S(X,X)=\n S(X, JX)+S\circ\n J(X,X).\tag 2.8$$
Thus
$$(S-\lb I)(\n J(X,X))=J(\n S(X,X))=-\frac12J\n \lb=0.\tag 2.9$$
From (2.4) we get $\n J(X,X)_{\mu}=0$ where $Y_{\mu}$ denotes the
$\2$-component of a vector  with respect to the decomposition
$TM=\1\oplus\2$. It is clear in view of (2.4) that $\n
J(X,X)_{\mu}=-g(X,X)(J\n\ln f)_{\mu}$ and thus $(\n f)_{\mu}=0$.
Thus $\xi\in\G(\1)$.

Set $TX=\n_X\xi$. Now, since $\1$ is totally geodesic it is clear
that $T\1\subset\1$ and $T\2\subset \2$.  Let $T^1=\n^1\xi$. Then
also $T^1\1\subset\1$ and $T^1\2\subset \2$. In fact,
$\n^1_X\xi=TX+d\ln f(X)\xi -g(X,\xi)\n\ln f$.

Now we shall prove that $JTX=\phi X$ for all $X\in \2$ for a
certain function $\phi$. Note that
$$\n S(X,\xi)=(\lb I-S)TX=(\lb-\mu)TX\tag 2.10$$
if $X\in \2$, and consequently
$$(\lb-\mu)g(TX,JX)=g(\n S(X,\xi),JX)=g(\n S(X,JX),\xi).\tag 2.11$$
Note that in view of (2.4) and (1.5) the following equality holds:
$$\gather \n S(X,JX)=(\mu I-S)(\n J(X,X)+J(\n_XX))\tag 2.12 \\
=(\mu-\lb)(- g(X,X)J(\n \ln f)+g(X,X)\frac
1{2(\lb-\mu)}J(\n\mu).\endgather $$ It follows that
$$g(TX,JX)=-g(JTX,X)=g(X,X)\bigg(-d\ln f(J\xi)+\frac
1{2(\lb-\mu)}d\mu(J\xi)\bigg).\tag 2.13$$ Thus for $X\in\2$,
$$JTX=\bigg(d\ln f(J\xi)-\frac 1{2(\lb-\mu)}d\mu(J\xi)\bigg)X.\tag 2.14$$
Let us recall the definition of a special K\"ahler-Ricci potential
([D-M-1], [D-M-2]).
 \medskip {\it Definition.}  A nonconstant
function $\tau\in C^{\infty}(M)$, where $\m$ is a K\"ahler
manifold, is called a special K\"ahler-Ricci potential if the
field $X= J(\n \tau)$ is a Killing vector field and at every point
with $d\tau\ne 0$ all nonzero tangent vectors orthogonal to the
fields $X,JX$ are  eigenvectors of both $\n d\tau$ and the Ricci
tensor $\rho$ of $\m$.

We now show  that the function $f$ is a special K\"ahler-Ricci
potential on the K\"ahler manifold $(M,g_1,J)$. In fact, for
$X\in\2$ we have
$$\gather \n^1_X\n^1f=-J(\n_X^1\xi)\\=-J(\n_X\xi+d\ln f(X)\xi+d\ln
f(\xi)X-g(X,\xi)\n\ln f)\\=-J(\n_X\xi)=-JTX=-\bigg(d\ln
f(J\xi)-\frac 1{2(\lb-\mu)}d\mu(J\xi)\bigg)X\endgather$$ and it is
clear that $\2$ is an eigendistribution of $\n^1df$. Since $\2$ is
also an eigendistribution of $\rho$ it follows from (2.6) that
$\2$ is also an eigendistribution of the Ricci tensor $\rho_1$ of
$(M,g_1,J)$. Thus on the whole of $U$ the function $f$ is a
special K\"ahler-Ricci potential. The fact that $f$ is a special
K\"ahler-Ricci potential on $V$ follows from the results of
Derdzi\'nski and Maschler ([D-M-1], [D-M-2]) note that $(V,g_1,J)$
is conformally Einstein.  Thus $f$ is a special K\"ahler-Ricci
potential on the open and dense subset $U\cup V$ of $M$ and
consequently is  a special K\"ahler-Ricci potential on the whole
of $M$.  If $\m$ is a K\"ahler manifold then it is a WBF manifold
(see [A-C-G-T]) and consequently it is extremal. We give a short
proof in this case for a completness and for the convenience of
the reader ( for another proof of this fact see [A-C-G]).  It
follows that $J(\n\tau)$ is a holomorphic Killing vector field
where $\tau$ is the scalar curvature of $(M,g,J)$. From our
assumptions it is clear that $\n\tau\in\G(\1)$. Thus
$$\n S(X,\n\tau)+(S-\lb I)(\n_X\tau)=0.\tag  2.15$$
On the other hand, note that every such K\"ahler metric satisfies
(see [J-1])
$$\gather \n_X\rho(Y,Z)=\frac 1{2\text{dim} M+4}(g(X,Y)Z\tau+g(X,Z)Y\tau+2g(Y,Z)X\tau\tag 2.16\\
-g(JX,Y)(JZ)\tau-g(JX,Z)(JY)\tau),\endgather$$ and consequently
$$\n S(X,Y)=\frac1{4(n+1)}(g(X,Y)\n\tau+Y\tau X-2Y
X\tau+g(JX,Y)J\n\tau-(JY\tau)JX)\tag 2.17$$ and $\n
S(X,\tau)=\frac
1{4n+4}(-X\tau\n\tau+|\n\tau|^2X+(JX\tau)J\n\tau)$.  Thus if
$X\in\2$ we have
$$(\mu-\lb)(\n_X\n\tau)_{\mu}=\frac1{4n+4}|\n\tau|^2X. $$
Since  $\1$ is totally geodesic and $\n\tau\in \1$ it follows that
$\n_X\n\tau\in\2$ if $X\in \2$ (note that
$g(\n_X\n\tau,Y)=g(\n_Y\n\tau,X)=H^{\tau}(X,Y)$). Consequently,
for $X\in\2$ we get
$$\n_X\n\tau=\frac1{(4n+4)(\mu-\lb)}|\n\tau|^2X, \tag 2.18$$
which means that $\tau$ is a special K\"ahler-Ricci potential.

From the results of Derdzi\'n\-ski and Maschler ([D-M-1], [D-M-2])
it follows that $(M,g,J)$ is in both cases a $\Bbb{CP}^1$-bundle
over a K\"ahler Einstein manifold $N$  or is $\Bbb{CP}^{n}$ in the
non-K\"ahler case. These $\Bbb{CP}^1$-bundles are
$\Bbb{P}(L\oplus\Cal O)$ where $L\rightarrow N$ is a complex line
bundle whose Euler class is proportional to the K\"ahler class of
$N$ and $\Cal O$ is the trivial line bundle.

\bigskip

\bigskip
{\bf 3. Construction of $\Cal{AC}^{\perp}$ Hermitian manifolds.}
In our construction we shall follow L. B\'erard Bergery (see
[Ber], [S]).
 Let $(N,h,J)$ be a compact K\"ahler
Einstein manifold, which is not Ricci flat and $\dim N=2m$,
$s\ge0,L>0,s\in \Bbb Q,L\in \Bbb R$, and $g:[0,L]\rightarrow \Bbb
R$ be a positive, smooth function on $[0,L]$  which is even at $0$
and $L$, i.e. there exists an $\e>0$ and  even, smooth functions
$g_1,g_2:(-\e,\e)\rightarrow \Bbb R$ such that $g(t)=g_1(t)$ for
$t\in[0,\e)$ and $g(t)=g_2(L-t)$ for $t\in(L-\e,L]$. Let
$f:(0,L)\rightarrow \Bbb{R}$ be positive on $(0,L)$, $f(0)=f(L)=0$
and let $f$ be odd at the points $0,L$. Let $P$ be a circle bundle
over $N$ classified by the integral cohomology class $\frac
s2c_1(N)\in H^2(N,\Bbb R)$ if $c_1(M)\ne 0$. Let $q$ be the unique
positive integer such that $c_1(N)=q\a$ where $\a\in H^2(N,\Bbb
R)$ is an indivisible integral class. Such a $q$ exists if $N$ is
simply connected or $\dim N=2$.  Note that every K\"ahler Einstein
manifold with positive scalar curvature is simply connected.  Then
$$s=\frac{2k}q; k\in\Bbb Z.$$ It is known that $q=n$ if
$N=\Bbb{CP}^{n-1}$ (see [Bes], p.273). Note that
$c_1(N)=\{\frac1{2\pi}\rho_N\}=\{\frac{\tau_N}{4m\pi}\0_N\}$ where
$\rho_N=\frac{\tau_N}{2m}\0_N$ is the Ricci form of $(N,h,J)$,
$\tau_N$ is the scalar curvature of $(N,h)$ and $\0_N$ is the
K\"ahler form of $(N,h,J)$. We can assume that $\tau_N=\pm 4m$. In
the case $c_1(N)=0$ we shall assume that $(N,h,J)$ is a Hodge
manifold, i.e. the cohomology class $\{\frac s{2\pi}\0_N\}$ is an
integral class. On the bundle $p:P\rightarrow N$ there exists a
connection form $\th$ such that $d\th=sp^*\0_N$ where
$p:P\rightarrow N$ is the bundle projection. Let us consider the
manifold $U_{s,f,g}=(0,L)\times P$ with the metric
$$g=dt^2+f(t)^2\th^2+g(t)^2p^*h.\tag 3.1$$
 It is known that the metric (3.1) extends to a metric
on the sphere bundle $M=P\times_{S^1}\Bbb{CP}^1$ if and only if a
function $g$ is positive and smooth on $[0,L]$, even at the points
$0,L$, the function $f$ is positive on $(0,L)$, smooth and odd at
$0,L$ and additionally $$f'(0)=1,\qquad f'(L)=-1\tag 3.2$$ Then
the metric (3.1) is bi-Hermitian.  We shall prove this in Section
4.

The metric (3.1) extends to a metric on $\Bbb{CP}^n$ if and only
if the function $g$ is positive and smooth on $[0,L)$, even at
$0$, odd at  $L$, the function $f$ is positive, smooth and odd at
 $0,L$ and additionally
$$f'(0)=1,\qquad f'(L)=-1,\qquad g(L)=0,\qquad g'(L)=-1.\tag 3.3$$

{\bf 4.  Circle bundles.  } Let $(N,h,J)$ be a K\"ahler manifold
with integral class  $\{\frac1{\pi}\0_N\}$ and let $p:P\rightarrow
N$ be a circle bundle with a connection form $\th$ such that
$d\th=s\0_N$, where $s\in\Bbb Q$ (see [K]). Let us consider a
Riemannian metric $g$ on $P$ where
$$g=a^2\th\otimes\th+b^2p^*h\tag 4.1$$
where $a,b\in \Bbb R$. Let $\xi$ be a fundamental vector field of
the action of $S^1$ on $P$, i.e. $\th(\xi)=1, L_{\xi}g=0$. It
follows that $\xi\in\frak{iso}(P)$ and $a^2\th=g(\xi,.)$.
Consequently,
$$a^2d\th(X,Y)=2g(TX,Y)\tag 4.2$$
for every $X,Y\in\frak X(P)$ where $TX=\n_X\xi$. Note that
$g(\xi,\xi)=a^2$ is constant, hence $T\xi=0$. On the other hand,
$d\th(X,Y)=sp^*\0_N(X,Y)=sh(Jp(X),p(Y))$.  Note that there exists
a tensor field $\tilde J$ on $P$ such that $\tilde J\xi=0$ and
$\tilde J(X)=(JX_*)^H$ where $X=X_*^H$ is the horizontal lift of
$X_*$. Indeed, $L_{\xi}T=0$ and $T\xi=0$, hence $T$ is a
horizontal lift of the tensor $\tilde T$. Now $\tilde
J=\frac{2b^2}{sa^2}\tilde T$.  Since $T\xi=0$ we get $\n
T(X,\xi)+T^2X=0$ and $R(X,\xi)\xi=-T^2X$. Thus
$g(R(X,\xi)\xi,X)=\|TX\|^2$ and
$$\rho(\xi,\xi)=\|T\|^2=\frac{sa^4}{4b^4}2m.$$ Consequently,
$$\lb=\rho\bigg(\frac{\xi}a,\frac{\xi}a\bigg)=\frac1{a^2}\|T\|^2=\frac{s^2a^2}{4b^4}m.\tag 4.3$$

We shall compute the O'Neill tensor $A$ (see [ON]) of the
Riemannian submersion $p:(P,g)\rightarrow (N,b^2h)$. We have
$$A_EF=\Cal V(\n_{\Cal H E}\Cal H F)+\Cal H(\n_{\Cal H E}\Cal V
F).$$

Let us write $u=\Cal V(\n_{\Cal H E}\Cal H F)$ and $v=\Cal
H(\n_{\Cal H E}\Cal V F)$.  The vertical component of a field $E$
equals $\th(E)\xi$. If $X,Y\in\Cal H$ then
$$g\bigg(\n_XY,\frac1a\xi\bigg)=\frac1a
(Xg(Y,\xi)-g(Y,\n_X\xi))=-\frac1ag(TX,Y)=\frac1ag(X,TY).\tag 4.4$$
Hence $u=\frac1ag(E-\th(E)\xi$ and
$T(F-\th(F)\xi)\frac{\xi}a=\frac1{a^2}g(E,TF)\xi.$ Note that $\Cal
H(\n_Xf\xi)=f\Cal H(\n_X\xi)=fTX$, hence $$v=\Cal H(\n_{\Cal H
E}\Cal V F)=\th(F)T(E)=\frac1{a^2}g(\xi,F)TE.$$ Consequently,
$$A_EF=\frac1{a^2}(g(E,TF)\xi+g(\xi,F)TE).\tag 4.5$$
If $U,V\in \Cal H$ then
$$\|A_UV\|^2=\frac1{a^2}g(E,TF)^2=\frac{s^2a^2}{4b^4}g(E,\tilde J
F)^2.$$ If $E$ is horizontal and $F$ is vertical then
$$A_EF=\frac1{a^2}g(\xi,F)TE.\tag 4.6$$
Hence $A_E\xi=TE$ and $\|A_E\xi\|^2=\|TE\|^2=\frac{s^2a^4}{4b^4}$.
It follows that
$$K(P_{E\xi})=\frac{s^2a^2}{4b^4},$$
where $K(P_{EF})$ denotes the sectional curvature of the plane
generated by the vectors $E,F$.  If $E,F\in\Cal H$ then
$$K(P_{EF})=K_*(P_{E_*F_*})-\frac{3g(E,TF)^2}{a^2\|E\w F\|^2},$$
where $E_*$ means the  projection of $E$ on $M$, i.e. $E_*=p(E)$.
Thus
$$K(P_{EF})=\frac1{b^2}K_0(P_{E_*F_*})-\frac{3s^2a^2g(E,\tilde JF)^2}{4b^4\|E\w F\|^2},\tag 4.7$$
where $K_0$ means the  sectional curvature of the metric $h$ on
$N$. Applying this we get for any $E\in\Cal H$ the formula for the
Ricci tensor $\rho$ of $(M,g)$:
$$\rho(E,E)=\frac1{b^2}\rho_0(bE_*,bE_*)-\frac{3s^2a^2}{4b^4}+\frac{s^2a^2}{4b^4},$$
where $\rho_0$ is the Ricci tensor of $(M,h)$. Hence
$$\mu=\frac{\mu_0}{b^2}-\frac{s^2a^2}{2b^4},$$
where $\rho_0=\mu_0g_0$. Now we shall find a formula for
$R(X,\xi)Y$ where $X,Y\in \Cal H$. We have $R(X,\xi)Y=\n T(X,Y)$
and $$\gather \n T(X,Y)=\n_X(T(Y))-T(\n_XY)=\n^*_{X_*}(\tilde
TY^*)+\frac12\Cal V[X,TY]\tag 4.8\\-(\tilde
T(\n^*_{X_*}Y_*))^*=\frac12\Cal
V[X,TY]=-\frac12sp^*\0_N(X,TY)\xi=-\frac{s^2a^2}{4b^2}h(X_*,Y_*)\xi.\endgather$$
Consequently, $R(X,Y,Z,\xi)=0$ for $X,Y,Z\in\Cal H$, and
$$R(X,\xi,Y,\xi)=-\frac{s^2a^4}{4b^2}h(X_*,Y_*).\tag 4.9$$

{\bf 5.   Riemannian submersion $  p:(0,L)\times P\rightarrow
(0,L)$.}  In this case the O'Neill tensor $A$ is $0$. We shall
compute the O'Neill tensor $T$ (see [ON]). Let $H=\frac d{dt}$ be
the horizontal vector field for this submersion and $\De$ be the
distribution spanned by the vector fields $H,\xi$. If $U,V\in \Cal
V$ and $g(U,V)=0$ then $T(U,V)=0$. We also have for $U\in \Cal V$
 with $g(U,\xi)=0$ and $U=U_*^*$ where $h(U_*,U_*)=1$,
$$T(U,U)=-gg'H.\tag 5.1$$
In fact, $2g(\n_UV,H)=-Hg(U,V)=-2gg'h(U_*,V_*)$ if $U=V$ and $0$
if $g(U,V)=0$.
 We also have
$$T(\xi,\xi)=-ff'H.\tag 5.2$$
Now we shall prove that the almost complex structure $J= J_{\e}$
defined by
$$JH=\e\frac1f\xi,JX=(J_*X_*)^*\qquad\text{   for  } X=(X_*)^*\in\E=\De^{\perp}$$
where $X_*\in TN,\e\in\{-1,1\}$, is  Hermitian   with respect to
the metric $g$.  For horizontal lifts $X,Y\in\frak
X(P)\subset\frak X((0, L)\times P)$ of the fields $X_*,Y_*\in\frak
X(N)$ ( with respect to the submersion described in  Section 4) we
have
$$\gather \n
J(Y,X)=\n_Y(JX)-J(\n_YX)=\n^*_{Y_*}(J_*(X_*))^*-\frac12d\th(Y,JX)\xi\\
+T(Y,JX)-J\bigg(\n^*_{Y_*}(X_*)^*-\frac12d\th(Y,X)\xi+T(Y,X)\bigg)\\
=-\frac12sh(JY,JX)fJH-gg'h(Y,JX)H-\frac12sh(JY,X)fH+h(X,Y)gg'JH\\
=h(X,Y)\bigg(\e gg'-\frac12 sf\bigg)JH+h(JY,X)\bigg(gg'-\e \frac12
sf\bigg)H.
\endgather$$
Hence $$\n J(JX,JY)=\n J(X,Y).\tag 5.4$$
 Since the distribution $\De$ is
totally geodesic and two-dimensional  it is clear that $\n
J(X,Y)=0$ if $X,Y\in\G(\De)$.  Now we shall show that
$$\n J(JX,JH)=\n J(X,H)\qquad\text {  for   } X\in \DE.$$
It is easy to show that $\n_XH=\n_HX=\frac{g'}gX$ and
$$\n_X(JH)=\e\n_X\bigg(\frac1f\xi\bigg)=\e\frac1fT(X)=\e\frac{sf}{2g^2}JX.\tag 5.5$$
On the other hand,
$$\n_X(JH)=\n J(X,H)+J(\n_XH)=\n J(X,H)+\frac {g'}gJ(X).\tag 5.6$$
Thus $\n J(X,H)=(\e\frac{sf}{2g^2}-\frac{g'}g)JX$. On the other
hand,
$$\n_X(JJH)=-\frac{g'}gX=\n_XJ(JH)+J(\n_XJH)=
\n J(X,JH)-\e\frac{sf}{2g^2}X.$$ Consequently,
$$\n J(X,JH)=\bigg(\e\frac{sf}{2g^2}-\frac{g'}g\bigg)X.$$
It  follows that $$\n J(JX,JH)=\n J(X,H).\tag 5.7$$

 Similarly
$$\n_H(JX)=\n_{JX}H=\frac{g'}gJX=\n J(H,X)+J(\n_HX)=\n
J(H,X)+\frac{g'}gJ\tag 5.8$$ and $\n J(H,X)=0$.  Let us recall
that $\n_X(\frac 1f\xi)=\frac{sf}{2g^2}JX$.  Thus we have
$$\n
J(JH,X)=\n_{JH}JX-J(\n_{JH}X)=\e\bigg(\n_{JX}\bigg(\frac1f\xi\bigg)-J\bigg(\n_X\bigg(\frac1f\xi\bigg)\bigg)=0.\tag
5.9
$$ Thus $$\n J(JH,JX)=\n J(H,X)\tag 5.10$$ and consequently  $$\n J(Y,Z)=\n
J(JY,JZ)$$ for all $Y,Z\in\frak X(M)$, which means that
$(M,g,J_{\e})$ is a Hermitian manifold.

 Let $U,V,W\in\Cal V$ and
$g(U,\xi)=g(V,\xi)=g(W,\xi)=0$. Then
$$R(U,V,\xi,W)=\hat
R(U,V,\xi,W)-g(T(U,\xi),T(V,W))+g(T(V,\xi),T(U,W))=0.$$
Analogously
$$R(U,V,\xi,H)=g(\n_VT(U,\xi),H)-g(\n_UT(V,\xi),H).\tag 5.11$$
Note that $T(U,\xi)=0$. Hence
$$0=g(\n_VT(U,\xi),H)+g(T(\n_VU,\xi),H)+g(T(U,\xi),\n_VH).$$
Consequently,
$$g(\n_VT(U,\xi),H)-g(\n_UT(V,\xi),H)=g(T([U,V],\xi),H)=0.\tag
5.12$$ From the  O'Neill formulae it also follows  that
$$R(JH,U,V,JH)=0$$ if $g(U,V)=0$, and
$$R(JH,U,U,JH)=\frac{s^2f^2}{4g^4}-\frac{f'g'}{fg}\tag 5.13$$ for
a unit vector field $U$ as above. Note also that the distribution
$\De$ spanned by the vector fields $\xi,H$ is totally geodesic.
Consequently, if $X,Y,Z\in\G(\De)$ and $V$ is as above then
$$R(X,Y,Z,V)=0.\tag 5.14$$

Note also that

$$ R(H,\frac1f\xi,\frac1f\xi,H)=-\frac {f''}f\tag 5.15$$
and
$$ R(H,X,X,H)=-\frac{g''}g\tag 5.16$$
for $X\in\G(\DE)$ and $||X||=1$.
\bigskip
{\bf 6.  Eigenvalues of the Ricci tensor.}  Let us assume that
$(N,h)$ is a $2(n-1)$-dimensional  K\"ahler-Einstein manifold of
scalar curvature $4(n-1)\e$ where $\e\in\{-1,0,1\}$.  Using the
results in Sections 3 and 4  we obtain the following formulae for
the eigenvalues of the Ricci tensor $\rho$ of
$(U_{s,f,g},g_{f,g})$:
$$\gather
\lb_0=-2(n-1)\frac{g''}g-\frac{f''}f,\tag 6.1\\
\lb_1=-\frac{f''}f+2(n-1)\bigg(\frac{s^2f^2}{4g^4}-\frac{f'g'}{fg}\bigg),\\
\lb_2=-\frac{g''}g+\bigg(\frac{s^2f^2}{4g^4}-\frac{f'g'}{fg}\bigg)+\frac{2\e}{g^2}-\frac{3s^2f^2}{4g^4}-(2n-3)\frac{(g')^2}{g^2}.\endgather$$

Since $\lb_0=\lb_1$ we obtain
$$\frac{g''}g+\frac{s^2f^2}{4g^4}-\frac{f'g'}{fg}=0.\tag 6.2$$
Hence $$f=\pm \frac {gg'}{\sqrt{\frac{s^2}4+Ag^2}},\tag 6.3$$
where $A\in \Bbb R$.  Using a homothety we can assume that $A\in
\{-1,0,1\}$. Under this condition the Ricci tensor $\rho$ has two
eigenvalues $\lb=-2(n-1)\frac{g''}g-\frac{f''}f$ and
$\mu=-\frac{g''}g-\frac{f'g'}{fg}+\frac{2\e}{g^2}-\frac{s^2f^2}{2g^4}-(2n-3)\frac{(g')^2}{g^2}$.
We shall discuss separately the conditions $A=0$ and $|A|=1$. Let
us write $h^2=\frac{s^2}4+Ag^2$. Then $g=\sqrt{|(\frac
s2)^2-h^2|}$ and $\operatorname{im} h\subset(-\frac s2,\frac s2)$
if $A=-1$ and $\operatorname{im} h\subset (\frac s2,\infty)$ if
$A=1$. Note that $h'\ne 0$. Let us assume that $h$ is an
increasing function. Then
$$f=\pm \frac{gg'}{|h|}=\pm \frac{hh'A}{|h|}= h'.$$
We also have
$$g''=\frac{A((h')^2+hh'')}g-\frac{hh'Ag'}{g^2}=\frac{A((h')^2+hh'')}g-\frac{h^2(h')^2}{g^3}.$$
Thus
$$\frac{g''}g=\frac{A((h')^2+hh'')}{g^2}-\frac{h^2(h')^2}{g^4}=-\frac{(h')^2+hh''}{r^2-h^2}-\frac{h^2(h')^2}{(r^2-h^2)^2}.$$
 Since $h'>0$ it follows that $t$ is a smooth function of $h$,
i.e. $t=t(h)$ and $\frac{dt}{dh}=\frac1{h'}$.  Let $z$ be a
function such that $h'=\sqrt{z(h)}$, i.e. $z(h)=h'(t(h))^2$. Then
$$f'=h''=\frac12z'(h).$$
Consequently,
$$\frac{g''}g=-\frac{z(h)+\frac12hz'(h)}{r^2-h^2}-\frac{h^2z(h)}{(r^2-h^2)^2}.$$
Analogously
$$\bigg(\frac{g'}g\bigg)^2=\bigg(\frac{hh'A}{g^2}\bigg)^2=\frac{h^2z(h)}{(r^2-h^2)^2},$$
and
$$\frac{f'g'}{fg}=\frac{\frac12z'(h)hh'A}{h'g^2}=\frac{\frac12z'(h)hA}{-A(r^2-h^2)}=-\frac{\frac12z'(h)h}{r^2-h^2}.$$
We also have $$\frac{s^2f^2}{2g^4}=2r^2\frac{z(h)}{(r^2-h^2)^2}.$$
Since
$$\mu=-\frac{g''}g-\frac{f'g'}{fg}+\frac{2\e}{g^2}-\frac{s^2f^2}{2g^4}-(2n-3)\frac{(g')^2}{g^2},$$
we get
$$\mu=\frac{hz'(h)}{r^2-h^2}-\frac{z(h)}{(r^2-h^2)^2}(r^2+(2n-3)h^2)-\frac{2\e A}{r^2-h^2}.$$

We shall use the following:
\medskip
{\bf Proposition 2. } {\it Let $\De$ be a distribution on
$U_{s,f,g}$  spanned by the fields $\{\th^{\sharp},H\}$. Then
$\De$ is a totally geodesic foliation with respect to the metric
$g_{f,g}$. The distribution $\DE$ is umbilical with the mean
curvature normal $\xi=-\n\ln g$. Let $\lb, \mu$ be the eigenvalues
of the Ricci tensor $S$ of $g_{f,g}$ corresponding to the
eigendistributions $\De,\De^{\perp}$ respectively.  Then the
following conditions are equivalent:

(a) There exists $E\in \Bbb R$ such that $\lb-\mu =Eg^2$,

(b) There exist $C,D\in \Bbb R$ such that $\mu =Cg^2+D$,

(c) $\lb-2\mu$ is constant,

(d) $(U_{s,f,g},g_{f,g})$ is a Hermitian Gray manifold.}
\medskip
{\it Proof.  } The first part of the assertion  is a consequence
of [J-2]. Note that $\n\lb=H\lb H,\n\mu=H\mu H$. Consequently,
$\operatorname{tr}_g\n S=\frac12\n\tau=(H\lb+(n-1)H\mu)H$. On the
other hand, one can easily check that $\operatorname{tr}_g\n
S=2(n-1)(\mu-\lb)\xi+H\lb H$. Thus

$$\frac{\n\mu}{2(\lb-\mu)}=\n\ln g.$$

Now we prove that (a) $\Rightarrow$ (b). If (a) holds then $\n
\mu=2Eg^2\frac{\n g}g=E\n g^2$. Thus $\n(\mu+Eg^2)=0$, which
implies (b).

(b)$\Rightarrow$ (a). We have
$$\frac{\n g}g=\frac{\n\mu}{2(\lb-\mu)}=\frac{Cg\n g}{\lb-\mu},$$
and consequently $\n g(\frac{Cg^2-(\lb-\mu)}{g(\lb-\mu)})=0$,
which is equivalent to (a).

(a)$\Rightarrow $(c). We have $\lb-\mu=Eg^2$ and consequently  $\n
\mu=2Eg\n g=E\n g^2$. Thus $\n \lb=\n (\mu+Eg^2)=2E\n g^2$ and
$\n\lb-2\n\mu=0$, which gives (c).

(c)$\Rightarrow$ (a). If $\n\lb=2\n\mu$ then
$\n\lb=4(\lb-\mu)\frac{\n g}g$. Consequently,
$\n\lb-\n\mu=2(\lb-\mu)\frac{\n g}g$ and $\n\ln|\lb-\mu|=2\frac{\n
g}g=2\n\ln g$, which means that $\n\ln|\lb-\mu|g^{-2}=0$.  It
follows that $\ln\frac{|\lb-\mu|}{g^2}=C$ for some $C\in\Bbb R$,
which is equivalent to  (a).

(d)$\Leftrightarrow$(c).  This equivalence follows from
Proposition 1. $\k$
\medskip
{\it Remark.}  From the condition (c) in Proposition 2 it is clear
that the eigenvalues  $\lb,\mu$ of the Ricci tensor are everywhere
different in the case of a $\Bbb{CP}^1$-bundle and coincide at an
exactly one point if $M=\Bbb{CP}^n$.
\bigskip

{\bf 7. Solutions of the linear ODE related to Gray manifolds.}
Using  Prop.  2 we see that the equation characterizing our Gray
manifold is
$$\frac{hz'(h)}{r^2-h^2}-\frac{z(h)}{(r^2-h^2)^2}(r^2+(2n-3)h^2)=\frac{2\e A}{r^2-h^2}+C(r^2-h^2)+D,$$
or
$$z'(h)-\frac{z(h)}{h(r^2-h^2)}(r^2+(2n-3)h^2)=\frac{2\e A}{h}+C\frac{(r^2-h^2)^2}h+D\frac{r^2-h^2}h,\tag 7.1$$
where $C,D$ are any real numbers. We are looking for the solutions
of (7.1) satisfying the following boundary conditions: if $A=-1$
then there exist real numbers $x,y$ such that $-r<x<y<r$ and if
$A=1$ then $r<x<y$ and in both cases
$z(x)=0,z'(x)=2,z(y)=0,z'(y)=-2$ and $z(t)>0$ if $t\in (x,y)$.

{\bf The case $s\ne 0,A\in\{-1,1\}$.}  Let us write $z
(t)=z_0(\frac tr)$. Then it is easy to check that the function
$z_0$ satisfies the equation
$$z_0'(h)-\frac{z_0(h)}{h(1-h^2)}(1+(2n-3)h^2)=\frac{2\eta}{h}+r^4C\frac{(1-h^2)^2}h+r^2D\frac{1-h^2}h,\tag 7.2$$
where $\eta=\e A$, and the boundary conditions, if $A=-1$ then
there exist real numbers $x,y$ such that $-1<x<y<1$  and if $A=1$
then $1<x<y$ and in both cases
$z_0(x)=0,z_0'(x)=2r,z_0(y)=0,z_0'(y)=-2r$ and $z(t)>0$ if $t\in
(x,y)$.

We have $z_0(t)=F(t)\frac t{(1-t^2)^{n-1}}$ and the function $F$
satisfies the equation:
$$F'(t)=2\eta\frac{(1-t^2)^{n-1}}{t^2}+C\frac{(1-t^2)^{n+1}}{t^2}+D\frac{(1-t^2)^{n}}{t^2},\tag 7.3$$
where we denote $r^4C,r^2D$ again by $C,D$. Hence $$\gather
F(t)=-\frac{2\eta}t+2\eta\sum_{k=1}^{n-1}(-1)^k
C^k_{n-1}\frac{t^{2k-1}}{2k-1}\tag 7.4\\-\frac
Ct+C\sum_{k=1}^{n+1}(-1)^k C^k_{n+1}\frac{t^{2k-1}}{2k-1}-\frac
Dt+D\sum_{k=1}^{n}(-1)^k
C^k_{n}\frac{t^{2k-1}}{2k-1}+E,\endgather$$ where $E\in \Bbb R$
and $C^k_n=\frac{n!}{k!(n-k)!}$. From the boundary conditions we
obtain
$$\gather C=\frac{2(-r+\eta x-\eta y+rxy)}{(x^2-1)(-x+y+xy^2-y^3)},\tag 7.5\\
 D=\frac{2(\eta(x^3-x^2y+x(y^2-2)-y(y^2-2)+r(1+x^3y-x^2y^2+xy(y^2-2))}{(x^2-1)(y^2-1)(x-y)},\endgather$$
 if $x\ne -y$. Consequently, the condition on $x,y$ such that $x\ne -y$ is as follows:

$$\gather-\frac{2\eta}x+2\eta\sum_{k=1}^{n-1}(-1)^k
C^k_{n-1}\frac{x^{2k-1}}{2k-1}-\frac Cx+C\sum_{k=1}^{n+1}(-1)^k
C^k_{n+1}\frac{x^{2k-1}}{2k-1}\tag 7.6\\-\frac
Dx+D\sum_{k=1}^{n}(-1)^k
C^k_{n}\frac{x^{2k-1}}{2k-1}+\frac{2\eta}y-2\eta\sum_{k=1}^{n-1}(-1)^k
C^k_{n-1}\frac{y^{2k-1}}{2k-1}\\+\frac Ct-C\sum_{k=1}^{n+1}(-1)^k
C^k_{n+1}\frac{y^{2k-1}}{2k-1}+\frac Dy-D\sum_{k=1}^{n}(-1)^k
C^k_{n}\frac{y^{2k-1}}{2k-1}=0,\endgather$$ where $C,D$ are given
by formulae (7.5) and
$$\gather-E=-\frac{2\eta}x+2\eta\sum_{k=1}^{n-1}(-1)^k
C^k_{n-1}\frac{x^{2k- 1}}{2k-1}-\frac Cx+C\sum_{k=1}^{n+1}(-1)^k
C^k_{n+1}\frac{x^{2k-1}}{2k-1}\\-\frac Dx+D\sum_{k=1}^{n}(-1)^k
C^k_{n}\frac{x^{2k-1}}{2k-1}.\endgather$$ We are looking for
$(x,y)$ lying on the  algebraic curve given by (7.6). We show that
the function given by a solution $(x,y)$ of these equations really
gives the $\AC$ Hermitian manifold at least if $x>0$ or $y<0$. Let
us assume that either $0<x<y<1$, $1<x<y$ or $-1<x<y<0$ and there
exist two points $p,q\in (x,y)$ such that $F'(p)=F'(q)=0$. Then
$$ \gather 0=2\eta+C(1-p^2)^{2}+D(1-p^2),\\
0=2\eta+C(1-q^2)^{2}+D(1-q^2).\endgather$$ Thus
$$C=C(p,q)=\frac{2\eta}{(p^2-1)(q^2-1)},\qquad D=D(p,q)=\frac{2\eta(-2+p^2+q^2)}{(p^2-1)(q^2-1)}.$$
It follows that for a given $p$ such that $F'(p)=0$ there  exists
at least one $q\in (x,y),q\ne p$ such that $F'(q)=0$ since the
function $I\ni q\mapsto C(p,q)\in\Bbb R$ where $I=(-1,0),I=(0,1)$
or $I=(1,\infty)$ is injective.

\medskip
If $x=-y$ then $E=0$. Note that if $E=0$ then $z_0$ is an even
function. Now we consider the case $A=-1$, $-1<x<0$ and $y=-x$.
Then
$$C=\frac{2rx-2\eta-D(1-x^2)}{(1-x^2)^2}.$$
We shall find $D$ such that  $z_0(x)=0$, which means $F(x)=0$. We
have

$$\gather-\frac{2\eta}x+2\eta\sum_{k=1}^{n-1}(-1)^k
C^k_{n-1}\frac{x^{2k-1}}{2k-1}-\frac{2rx-2\eta-D(1-x^2)}{(1-x^2)^2x}\\
+\frac{2rx-2\eta-D(1-x^2)}{(1-x^2)^2}\sum_{k=1}^{n+1}(-1)^k
C^k_{n+1}\frac{x^{2k-1}}{2k-1}\\-\frac Dx+D\sum_{k=1}^{n}(-1)^k
C^k_{n}\frac{x^{2k-1}}{2k-1}=0.\endgather$$

Hence $$\gather
D\bigg(\frac1{1-x^2}-\frac1{1-x^2}\sum_{k=1}^{n+1}(-1)^k
C^k_{n+1}\frac{x^{2k-1}}{2k-1}-1+\sum_{k=1}^{n}(-1)^k
C^k_{n}\frac{x^{2k-1}}{2k-1}\bigg)\\
=2\eta-2\eta\sum_{k=1}^{n-1}(-1)^k
C^k_{n-1}\frac{x^{2k-1}}{2k-1}+\frac{2rx-2\eta}{(1-x^2)^2}-\frac{2rx-2\eta}{(1-x^2)^2}\sum_{k=1}^{n+1}(-1)^k
C^k_{n+1}\frac{x^{2k-1}}{2k-1}.\endgather$$

Thus

$$\gather
D\bigg(2x^2-C^2_{n+1}\frac{x^4}3+C^2_n\frac{x^4}3+nx^4+\phi_5(x)\bigg)\tag 7.7\\
=\frac1{1-x^2}(2rx+2\eta(n-3)x^2+2r(n+1)x^3+\psi_4(x))\endgather$$
and
$$D\bigg(1+\frac{n}3x^2+\phi_3(x)\bigg)=\frac1{1-x^2}\bigg(\frac
rx+\eta(n-3)+r(n+1)x+\psi_2(x)\bigg)\tag 7.8$$ and
$$D=\frac1{1-x^2}\bigg(\frac
rx+\eta(n-3)+2r\bigg(\frac56n+1\bigg)x+\psi_2(x)\bigg).\tag7.9$$
Hence
$$C=-\frac rx-2\eta -\eta(n-3)-2r\bigg(\frac56n+1\bigg)x+\phi_2(x),\tag 7.10$$
where $f_k(x)$ denotes an analytic function of the form
$f_k(x)=x^k(a_0+a_1x+a_2x^2+\ldots)$ convergent in a neighborhood
of $x=0$.  It follows that
$$\gather z_x(0)=-2\eta-C-D=-2\eta+\frac rx+2\eta
+\eta(n-3)+2r\bigg(\frac56n+1\bigg)x\tag 7.11\\-\phi_2(x)-\frac
rx-\eta(n-3)-2r\bigg(\frac56n+1\bigg)x-rx+\gamma_2(x))=-rx+\a_2(x).\endgather$$
Consequently, $z_x(0)>0$ for  sufficiently small $|x|$ where $x\in
(-1,0)$. Now we shall show that the function $z_x$ satisfies the
condition $z_x(t)>0$ for $t\in (x,-x)$. To this end we shall look
at $F_x'(t)=\frac{(1-t^2)^{n-1}}{t^2}(2\eta+C(1-t^2)+D(1-t^2))$
where $C,D$ are given by $(7.9),(7.10)$. We shall prove that for
small $x$ the function $F_x'(t)$ has no zeros in the interval
$(x,0)$. Since $F(x)=0,F'(x)\frac x{(1-x^2)^{n-1}}=2r>0$ and
$\lim_{t\rightarrow 0^-}F(t)=-\infty$ it follows that if $F(t)$
had a zero in $(x,0)$ then $F(t)$ would have at least two zeros in
$(x,0)$.

Note that $F'(t)=0$ if and only if (we shall assume that $C\ne 0$)
$$C\bigg(\a^2+\frac DC\a+\frac{2\eta}C\bigg)=0,\tag 7.12$$ where $\a=1-t^2$.
Note that
$$\frac DC=-\frac{1+\eta(n-3)\frac xr+\psi_2(x)}{1+\eta(n-1)\frac xr+\phi_2(x)}.$$
Consequently,
$$\frac
DC=-1+\frac{2\eta}rx+\g_2(x).$$ Analogously
$8\frac{\eta}C=8\eta\frac xr+\beta_2(x)$.  It follows that
$$\a=\frac12\bigg(1-\frac{2\eta}rx-\g_2(x)\pm\sqrt{1+\frac{4\eta}rx+o_2(x)}  \bigg)$$
are the solutions of (7.12). Hence one of the roots is
$\a=1-t^2=-4\frac{\eta}rx+k_2(x)$ and then
$t^2=1+4\frac{\eta}rx-k_2(x)$. For a small $x$ this $|t|$ is close
to $1$. Thus $F'(t)$ can have at most one zero in $(-x,0)$. It
follows that $F(t)$ is different from $0$ in $(x,0)$ and
consequently $z_x(t)$ is positive on the whole of $(x,-x)$ for
small $x\in (-1,0)$. This means that in the case of $|A|=1$ we
have infinitely many examples of compact  bi-Hermitian Gray
manifolds of any dimension $2n$ corresponding to a given compact
K\"ahler Einstein manifold $(N,h)$ and a rational number $s$
described in Section 3. They are  holomorphic $\Bbb{CP}^1$-bundles
over $N$.

Now we consider the case where $M=\Bbb{CP}^n$. Then
$N=\Bbb{CP}^{n-1}$, $s=\frac2n,k=1$. One of the zeros of $z_0$ is
$1$ and let $y$ be the other zero such that $z_0>0$ on the
interval $(y,1)$ or $(1,y)$. Note that $y<1$ if $\eta=A=-1$, and
$y>1$ if $A=1$. It is easily seen, that $z'_0(1)=\frac{2A}n=A s$.
Since $z_0(y)=0$ and $z'_0(y)=-sA=-\frac{2A}n$ we get from (7.2):
$$2A+C(1-y^2)^2+D(1-y^2)=-\frac{2Ay}n.\tag 7.13$$

Let us assume that $z'_0(x)=0$ for $x\in (y,1)$ (or $x\in (1,y)$
respectively). The functions $z_0(t),F(t)$ have the same set of
zeros in the interval $(0,+\infty)$. We show that the only zeros
of $F(t)$ in  $[y,1]$ (or $[1,y]$) are $1,y$. Since $F(y)=F(1)=0$
it is enough  to show that $F'(t)$ has at most one zero in $(y,1)$
($(1,y)$ respectively).  Let $F'(x)=0$. Then
$$-2A=C(1-x^2)^2+D(1-x^2),\tag 7.14$$
and

$$\gather
C=C(x,y)=\frac{2A(nx^2-y+x^2y-ny^2)}{n(1-x^2)(x^2-y^2)(1-y^2)},\tag 7.15\\
D=D(x,y)=\frac{2A(-2nx^2+nx^4+y-2x^2y+2ny^2-ny^4)}{n(y^2-1)(x^2-y^2)(x^2-1)}.\endgather$$

The derivative of $C=C(x,y)$ as a function of $x$ is $$\frac{\p
C}{\p
x}=\frac{4Ax((x^2-1)^2y+n(x^2-y^2)^2)}{n(x^2-1)^2(x^2-y^2)^2(1-y^2)},$$
and is always less than $0$ if $A=1$ and greater than $0$ if
$A=-1$ for $y>0$. It follows that the function $(y,1)\ni x\mapsto
C(x,y)\in \Bbb R$ is injective. It follows that $F'(t)$ can have
at most one zero in $(y,1)$. Consequently, the function $z_0$ is
positive on $(y,1)$ (or $(1,y)$) if $z_0$ is a solution of $(7.2)$
satisfying the boundary conditions   $
z_0(y)=0,z_0'(y)=-\frac{2A}n,z_0(1)=0$. Consequently, we obtain
two families, $(\Bbb{CP}^n,g_x),x\in (0,1)\cup(1,\infty)$, of
Hermitian $\AC$ metrics on $\Bbb{CP}^n$. Note that the solutions
exist also for $x<0$ (compare for example the case $n=2$ [J-3]).

\bigskip
{\bf The case $A=0$.}  This case can also be deduced from
[A-C-G-T]. For a convenience of the reader we give a  shorter
proof here.  Now we shall consider the case $A=0$, which
corresponds to the K\"ahler manifold $(M,g_{f,g})$. In this case
we have $f=\frac{2gg'}s$ and
$$\mu=-\frac{g''}g+\frac{2\e}{g^2}-2n\frac{(g')^2}{g^2},$$
and we obtain the  equation
$$z'(g)-2n\frac{z(g)}{g}=\frac{2\e }{g}-Cg^3-Dg,\tag
7.16$$

It follows that

$$z(g)=\frac{\e}n-C\frac{g^4}{2n+4}-D\frac{g^2}{2n+2}+\frac
E{g^{2n}}.$$ We are looking for the solution $z(t)$ satisfying for
certain points $0<x<y$ the initial conditions:

$$z(x)=z(y)=0,\qquad xz'(x)=s,\qquad yz'(y)=-s,\tag 7.17$$ and such that
$z(t)>0$ for $t\in(x,y)$. Hence

$$\gather C=\frac{-2\e+s-Dx^2}{x^4}, E=-\frac{x^{2n}(4\e+4\e
n+ns+n^2s-Dnx^2)}{2n(1+n)(2+n)},\\
 C=\frac{-2\e-s+Dy^2}{y^4}, E=\frac{y^{2n}(-4\e-4\e
n+ns+n^2s+Dny^2)}{2n(1+n)(2+n)}.\endgather$$

It follows that
$$\gather  D=\frac{2\e x^4+sx^4-2\e y^4+sy^4}{(x^2-y^2)x^2y^2}\\
D=\frac{(1+n)(4\e x^{2n}+nsx^{2n}-4\e
y^{2n}+nsy^{2n})}{n(x^{2+2n}-y^{2+2n})}.\endgather$$  Hence there
exists a function $z$ satisfying $(7.17)$ if and only if
$$\gather n(s+2\e)x^{2n+6}-n(s-2\e
)y^{2n+6}-[(6n+4)\e+n^2s]x^{2n+2}y^4\tag 7.18
\\-[(6n+4)\e-n^2s]y^{2n+2}x^4
+[(4\e(n+1)+n(n+1)s)]x^{2n+4}y^2\\+[(4\e(n+1)-n(n+1)s)]y^{2n+4}x^2=0\endgather$$
for certain $0<x<y$. This  equation is equivalent to

$$\gather -n(s-2\e)c^{2n+6}+[(4\e(n+1)-n(n+1)s)]c^{2n+4}-[(6n+4)\e-n^2s]c^{2n+2}\tag 7.19\\
-[(6n+4)\e+n^2s]c^4+[4\e(n+1)+n(n+1)s]c^2+n(s+2\e)=0,\endgather$$
where $c=\frac yx>1$. Let
$$\gather\phi(c)= -n(s-2\e)c^{2n+6}+[(4\e(n+1)-n(n+1)s)]c^{2n+4}-[(6n+4)\e-n^2s]c^{2n+2}\\
-[(6n+4)\e+n^2s]c^4+[4\e(n+1)+n(n+1)s]c^2+n(s+2\e).\endgather$$

Then $\phi(1)=0$ and $\phi'(1)=-8ns(n+1)<0$. It follows that if
$s-2\e<0$ then there exists a root $c_0\in (1,\infty)$ of equation
(7.19). Note that $\phi(c)=\e\Phi(c)+s\Psi(c)$, where $\Phi(c)=2n
c^{2n+6}+4(n+1)(c^{2n+4}+c^2)-(6n+4)(c^4+c^2)+2n$ and
$\Psi(c)=-n(c^{2n+6}+(n+1)(c^{2n+4}-c^2)-n(c^{2n+2}-c^4)-n)$. Note
that for $c>1$, $\Phi(c)>0$ and $\Psi(c)<0$. This follows from the
inequalities $(n+1)(c^{2n+4}-c^2)>n(c^{2n+2}-c^4)$ and
$2n(c^{2n+6}+1)+4(n+1)(c^{2n+4}+c^2)>(6n+4)(c^4+c^2)$ valid for
$c>1$. Hence for $\e\in\{-1,0\}$ we get $\phi(c)<0$ for $c>1$. On
the other hand, it is easy to see that $\Phi(c)+2\Psi(c)<0$ if
$c>1$. It follows that the polynomial $\phi$ has a real root
$c_0>1$ if and only if $\e=1$ and $s<2$.  We shall show that the
function
$$z(t)=\frac{\e}n-C\frac{t^4}{2n+4}-D\frac{t^2}{2n+2}+\frac
E{t^{2n}},$$ is positive on $(x,c_0x)$ for any $x>0$. Note that
the derivative of $z(t)$ can have at most two positive zeros and
consequently if $z(t)$ satisfies the conditions (7.17) then it is
positive for $t\in (x,y)$ where $y=c_0x$. Note also that for a
given  $c_0$ the solutions given by $(x,c_0x)$ give homothetic
metrics for all $x>0$. Since $s=\frac{2k}q<2$ it follows that
$k=1,2,...,q-1$ and for these values of $k$ there exists a
K\"ahler Gray bi-Hermitian metric on an appropriate
$\Bbb{CP}^1$-bundle over a compact K\"ahler Einstein manifold of
positive scalar curvature. Note that every such K\"ahler metric
satisfies the equation (see [J-1])
$$\gather \n_X\rho(Y,Z)=\frac 1{2\text{dim} M+4)}(g(X,Y)Z\tau+g(X,Z)Y\tau+2g(Y,Z)X\tau\\
-g(JX,Y)(JZ)\tau-g(JX,Z)(JY)\tau).\endgather$$
\medskip
{\bf The case $s=0$.} We shall assume that $g$ is not constant.
Then $f=g'$ and
$$\mu=-2\frac{g''}g-(2n-3)\frac{(g')^2}{g^2}+\frac{2\e}{g^2}.$$
 Thus we obtain the equation
$$z'(g)+(2n-3)\frac{z(g)}g=\frac{2\e}g-Cg^3-Dg.$$
It follows that
$$z(t)=\frac{2\e}{2n-3}-C\frac{t^4}{2n+1}-D\frac{t^2}{2n-1}+E\frac1{t^{2n-3}}.$$
Again the derivative of $z$ can have at most two positive roots.
Consequently, in the case $\e\ne 0$ there exists a nontrivial
solution satisfying the initial conditions
$$z(x)=z(y)=0,\qquad z'(x)=2,\qquad z'(y)=-2,\tag 7.20$$
for $0<x<y$ if
$$\gather
x=\e\frac{4(1-c)(c^{2n+1}-1)(2n-1)+2(-1+c-c^2+c^3)(c^{2n-1}-1)(2n+1)}{(2n-3)c((c^{2n+1}-1)(2n-1)+(1-c+c^2)(1-c^{2n-1})(2n+1))}
\endgather$$
is positive for some $c\in (1,+\infty)$. For $n=2,3$ it can happen
only if $\e=1$, and for  $n\in\{4,5,\ldots\}$, $x$ is positive for
some $c>1$ for both values of $\e\in\{-1,1\}$. Note that the
polynomial
$4(1-c)(c^{2n+1}-1)(2n-1)+2(-1+c-c^2+c^3)(c^{2n-1}-1)(2n+1)$ is
negative for $c>1$. So the sign of $\e x$ is opposite to the sign
of $Q(c)=c((c^{2n+1}-1)(2n-1)+(1-c+c^2)(1-c^{2n-1})(2n+1))$. If
$n\in\{2,3\}$ then $Q$ is strictly negative for $c>1$. For $n=2$
this is proved in [J-2],[J-3] and for $n=3$ we have
$Q(c)=-(c-1)^5c(2c^2+3c+2)$. Note also that $Q(1)=Q'(1)=Q''(1)=0$
and $Q^{(3)}(1)=2(1+2n)(2n^2-7n+3)$. Hence $Q(c)>0$ for $c>1$
close to $1$ if $n>3$. Since
$\lim_{c\rightarrow+\infty}Q(c)=-\infty$ it follows that $Q$
changes sign and consequently there exist solutions for both
$\e=1$ and  $\e=-1$. This means that if $(M,g,J)$ is a compact
K\"ahler Einstein manifold of dimension $\ge 6$ and of non-zero
scalar curvature then there exists a family of compact
bi-Hermitian Gray metrics on the manifold $\Bbb{CP}^1\times M$
with a warped product metric. If $\dim M\le 4$  then the
appropriate metric exists only if the scalar curvature of the
K\"ahler Einstein manifold $M$ is positive.

 In the
case $\e=0$ there exists a solution satisfying the boundary
conditions if the polynomial
$F(c)=(1+c)(1-c^{2n+1})(2n-1)+(1+c^3)(c^{2n-1}-1)(2n+1)$ has a
real root $c>1$. Then $n>3$ again. In fact, $F(1)=0$,
$F'(1)=F''(1)=0$ and $F^{(3)}(1)=(n-3)(1-4n^2)$. It follows that
for $n>3$ the polynomial $F$ takes negative values for $c>1$ close
to $1$. Since $\lim_{c\rightarrow+\infty}F(c)=+\infty$ it follows
that $F$ has a real root $c_0>1$. If $n=3$ then
$F(c)=(c-1)^5(2c^3+5c^2+5c+2)$ is obviously positive for $c>1$.
The case $n=2$ was considered in [J-3]. In that way we obtain Gray
bi-Hermitian metrics on the products $\Bbb{CP}^1\times M$ where
$M$ is a K\"ahler Ricci flat manifold of dimension $\ge 6$.

\medskip

\bigskip
\newpage
\centerline{\bf References.}
\par
\medskip
[A-C-G] V. Apostolov, D.M.J. Calderbank, P. Gauduchon Hamiltonian
2-forms in K\"ahler geometry, I general theory, J. Diff. Geom.
73,(2006), 359-412.
\par
\medskip
\cite{A-C-G-T} V. Apostolov, D.M.J. Calderbank, P. Gauduchon and
Ch. W. T\o nsen  {\it Hamiltonian 2-forms in K\"ahler geometry, IV
Weakly Bochner-flat K\"ahler manifolds}, Comm. in Analysis and
Geom. {\bf 16}, (2008), 91-126.

\par
\medskip
\cite{Bes}  A. L. Besse {\it Einstein manifolds}, Ergebnisse,
ser.3, vol. 10, Springer-Verlag, Berlin-Heidelberg-New York, 1987.
\par
\medskip
\cite{Ber}  L. B\'erard Bergery,{\it Sur de nouvelles vari\'et\'es
riemanniennes d'Einstein}, Pu\-bl. de l'Institute E. Cartan
(Nancy) {\bf 4},(1982), 1-60.
\par
\medskip
\cite{D}  A. Derdzi\'nski, {\it Classification of certain compact
Riemannian manifolds with harmonic curvature and non-parallel
Ricci tensor}, Math. Z. {\bf 172}  (1980), 273-280.
\par
\medskip
\cite{D-M-1} A. Derdzi\'nski, G. Maschler {\it Special
K\"ahler-Ricci potentials on compact K\" ahler manifolds},  J.
reine angew. Math. {\bf 593}, (2006), 73-116.
 \par
\medskip
\cite{D-M-2}  A. Derdzi\'nski, G.  Maschler {\it Local
classification of conformally-Einstein  K\"ahler metrics in higher
dimension},  Proc. London Math. Soc. (3) {\bf 87}, (2003), no. 3,
779-819.
 \par
\medskip
\cite{G} A. Gray {\it Einstein-like manifolds which are not
Einstein},  Geom. Dedicata {\bf 7}, (1978),  259-280.
\medskip
\cite{J-1} W. Jelonek,  {\it Compact K\"ahler surfaces with
harmonic anti-self-dual Weyl tensor}, Diff. Geom. and Appl. {\bf
16}, (2002),267-276.

\par
\medskip
\cite{J-2} W. Jelonek,  {\it Bi-Hermitian Gray surfaces }, Pacific
Math. J, {\bf 222},  no. 1, ( 2005),  57-68.
\par
\medskip
\cite{J-3} W. Jelonek {\it  Bi-Hermitian Gray surfaces II}, (to
appear in Diff. Geom and  Appl.)

 \par
\medskip
\cite{K} S. Kobayashi {\it Principal fibre bundles with the
1-dimensional toroidal group}, T\^ohoku Math.J. {\bf 8},(1956)
29-45.
\par
\medskip
\cite{K-N} S. Kobayashi and K. Nomizu {\it Foundations of
differential Geometry}, vol.2, Interscience, New York  1963
\par
\medskip
\cite{K-S} N. Koiso and Y. Sakane,{\it Non-homogeneous
K\"ahler-Einstein metrics on compact complex manifolds II}, Osaka
J. Math. {\bf 25} (1988),  933-959.
\par
\medskip
\cite{ON} B. O'Neill, {\it The fundamental equations of a
submersion}, Mich. Math. J. {\bf 13}, (1966),459-469.
\par
\medskip
\cite{S} P. Sentenac ,{\it Construction d'une m\'etriques
d'Einstein  sur la somme de deux projectifs complexes de dimension
2}, G\'eom\'etrie riemannienne en dimension 4 ( S\'eminaire Arthur
Besse 1978-1979) Cedic-Fernand Nathan, Paris (1981), pp. 292-307.
\par
\medskip
\cite{W-W} J. Wang and M. Wang Einstein metrics on $S^2$ bundles
Math. Ann. {\bf 310}, (1998), 497-526.
\par
\medskip

\par
\medskip
\newpage
Institute of Mathematics

Cracow University of Technology

Warszawska 24

31-155      Krak\'ow, POLAND.

 E-mail address: wjelon\@pk.edu.pl

\enddocument